\documentclass[12pt,a4paper,oneside,reqno]{article}

\usepackage[utf8]{inputenc}
\usepackage[english]{babel}
\usepackage[left=1in,right=1in,top=1in,bottom=1in]{geometry}
\usepackage{amsfonts,amssymb,amsmath,amsthm}
\usepackage{xcolor}
\usepackage[normalem]{ulem}
\usepackage{hyperref}
\hypersetup{colorlinks=true,linkcolor=blue,citecolor=blue,urlcolor=blue}
\usepackage{fancyhdr}
\usepackage{longtable}
\newtheorem{theorem}{Theorem}
\newtheorem{lemma}[theorem]{Lemma}
\newtheorem{corollary}[theorem]{Corollary}
\theoremstyle{remark}
\newtheorem{remark}[theorem]{Remark}
\newtheorem*{RemA}{Remark A (characteristic $3$)}
\newtheorem*{RemB}{Remark B (proof of Theorem 10)}
\newtheorem*{RemC}{Remark C (Theorem 15 in positive characteristic)}

\newcommand{\br}[2]{\{#1,#2\}}
\newcommand{\Js}{J(\Gamma,\{\,,\})}

\newcommand{\doi}[1]{\href{https://doi.org/#1}{doi:#1}}

\begin{document}

\noindent\fbox{\parbox{\dimexpr\textwidth-2\fboxsep-2\fboxrule\relax}{\small
\textbf{Corrected version verified in Lean.}  
This is a corrected version of the paper
V.~Zhelyabin and I.~Kaygorodov, \emph{On $\delta$-superderivations of simple superalgebras of
Jordan brackets}, St. Petersburg Math. J. 23 (2012), no.~4, 665--677, \doi{10.1090/S1061-0022-2012-01213-6}.}}

\bigskip

\begin{center}
{\Large On $\delta$-superderivations of simple superalgebras of Jordan brackets}

\bigskip

{\bf Victor Zhelyabin and Ivan Kaygorodov}
\end{center}

\bigskip

\noindent{\small\textbf{Abstract.} A complete description of $\delta$-derivations and
$\delta$-superderivations is given for simple unital superalgebras of Jordan brackets over a field of
characteristic different from~$2$ \textup{(}and, for the odd $\delta$-superderivations, different
from~$3$; see Remark~A for characteristic~$3$\textup{)} and for simple unital finite-dimensional Jordan
superalgebras over an algebraically closed field of characteristic $p\neq2$ \textup{(}with the
corrections of Remarks~A and~C for $p>2$\textup{)}. As a consequence, a criterion for simple unital
superalgebras of Jordan brackets to be special is obtained.}

\medskip

\noindent{\small\textbf{Keywords:} $\delta$-superderivation, $\delta$-derivation,
Jordan superalgebra, Kantor double.

}

\medskip

\bigskip

\setcounter{tocdepth}{1}
\tableofcontents

\section*{Introduction}
\addcontentsline{toc}{section}{Introduction}

The notion of a derivation of an algebra has been generalized by many mathematicians in various
directions. In particular, in~\cite{1} one can find the definition of a $\delta$-derivation of an
algebra. We recall that for fixed $\delta\in F$, a $\delta$-derivation of an $F$-algebra $A$ is a
linear map $\varphi$ satisfying the condition
\[
\varphi(xy)=\delta(\varphi(x)y+x\varphi(y))
\]
for any $x,y\in A$. In~\cite{1}, the $\frac12$-derivations were described for an arbitrary prime Lie
$F$-algebra $A$ ($\frac16\in F$) with a nondegenerate symmetric invariant bilinear form. Namely, it was
proved that a linear map $\varphi\colon A\to A$ is a $\frac12$-derivation if and only if
$\varphi\in\Gamma(A)$, where $\Gamma(A)$ is the centroid of~$A$. This implies that if $A$ is a central
simple Lie algebra over a field of characteristic $p\neq2,3$ with a nondegenerate invariant bilinear
form, then any $\frac12$-derivation $\varphi$ has the form $\varphi(x)=\alpha x$ for some $\alpha\in F$.

In~\cite{2} Filippov proved that any prime Lie $\Phi$-algebra admits no nonzero
$\delta$-derivations if $\delta\notin\{-1,0,\frac12,1\}$. Also in~\cite{2}, it
was shown that any prime Lie $\Phi$-algebra $A$ ($\frac16\in\Phi$) with a nonzero antiderivation is a
three-dimensional central simple algebra over the field of fractions of the center $Z_R(A)$ for its
algebra of right multiplications $R(A)$. In the same paper, a nontrivial $\frac12$-derivation was
constructed for the Witt algebra $W_1$, i.e., a $\frac12$-derivation that is not an element of the
centroid of~$W_1$.

The paper~\cite{3} contains a description for the $\delta$-derivations of the prime alternative and
non-Lie Maltsev $\Phi$-algebras with some restrictions on the operator ring~$\Phi$. It was shown that
the algebras in these classes have no nonzero $\delta$-derivations if $\delta\neq0,\frac12,1$. The
results of Filippov were partially generalized by Leger and Luks
in~\cite{9}. Those authors considered quasiderivations of Lie algebras, i.e., linear maps $f$ for which
there exists a linear map $f'$ related to $f$ by the formula $f'(xy)=f(x)y+xf(y)$. They proved that
the space of quasiderivations of a simple finite-dimensional Lie algebra $A$ of rank greater than~$1$
coincides with the direct sum of the space of derivations and the centroid of~$A$.

In~\cite{4}, a description of $\delta$-derivations was given for the simple finite-dimensional Jordan
superalgebras over an algebraically closed field of characteristic zero and for the semisimple
finite-dimensional Jordan algebras over an algebraically closed field of characteristic different
from~$2$. Later, in~\cite{5}, the $\delta$-derivations were described for the classical Lie
superalgebras. The paper~\cite{6} was devoted to a description of the $\delta$-derivations of Cartan
Lie superalgebras. In that paper, the $\delta$-superderivations were described for the simple
finite-dimensional Lie superalgebras. Also in~\cite{6}, the $\delta$-derivations of semisimple
finite-dimensional Jordan algebras and $\delta$-superderivations of Jordan superalgebras over an
algebraically closed field of characteristic zero were described. For algebras and superalgebras from
the papers~\cite{4,5,6}, the absence of nontrivial $\delta$-derivations and $\delta$-superderivations
was proved. Later, the results of~\cite{5} were generalized by Zusmanovich in~\cite{8}. He described the
$\delta$-derivations and $\delta$-superderivations of prime Lie superalgebras. He showed that a prime
Lie superalgebra has no nontrivial $\delta$-derivations and $\delta$-superderivations if
$\delta\notin\{-1,0,\frac12,1\}$. He proved that for a Lie superalgebra $A$ with zero center and a
nondegenerate supersymmetric invariant bilinear form satisfying the condition $A=[A,A]$, the space of
$\frac12$-derivations ($\frac12$-superderivations) coincides with the centroid (supercentroid) of the
superalgebra~$A$. Also, Zusmanovich gave an affirmative answer to the question of Filippov in~\cite{2}
concerning the existence of zero divisors in the ring of $\frac12$-derivations of a prime Lie algebra.
Subsequently, the $\delta$-superderivations of the generalized Kantor double constructed on a prime
nonunital associative algebra were considered in~\cite{7}.

In the present paper, we consider the $\delta$-derivations and $\delta$-superderivations in the case of
simple superalgebras of Jordan brackets. We prove the absence of nontrivial $\delta$-derivations and
$\delta$-superderivations of simple superalgebras of Jordan brackets, which are not superalgebras of
vector type. A description is given for the $\delta$-derivations and $\delta$-superderivations of
simple Jordan superalgebras of vector type. As a consequence, with the help of the classification of
all simple finite-dimensional Jordan superalgebras, as given in~\cite{13}, we obtain a description of
the $\delta$-derivations and $\delta$-superderivations of the simple unital finite-dimensional Jordan
superalgebras over an algebraically closed field of characteristic $p\neq2$ \textup{(}for $p>2$
see Remarks~A and~C\textup{)}.

\section{The main facts and definitions}

Let $F$ be a field of characteristic $p\neq2$. An algebra $A$ over the field $F$ is said to be Jordan
if it satisfies the identities
\[
xy=yx,\qquad(x^2y)x=x^2(yx).
\]
Let $G$ be a Grassmann algebra over $F$, given by generators $1,\xi_1,\dots,\xi_n,\dots$ and the
defining relations $\xi_i^2=0$, $\xi_i\xi_j=-\xi_j\xi_i$. The unity $1$ and the products
$\xi_{i_1}\xi_{i_2}\cdots\xi_{i_k}$, $i_1<i_2<\dots<i_k$, form a basis of the algebra $G$ over~$F$.
Denote by $G_0$ and $G_1$ the subspaces generated by products of even and odd length, respectively;
then $G$ is a direct sum of these subspaces: $G=G_0\oplus G_1$; moreover, we have the inclusions
$G_iG_j\subseteq G_{i+j\,(\mathrm{mod}\,2)}$, $i,j=0,1$. In other words, $G$ is a $\mathbb Z_2$-graded
algebra (or a superalgebra) over~$F$.

Now, let $A=A_0\oplus A_1$ be an arbitrary superalgebra over~$F$. Consider the tensor product
$G\otimes A$ of $F$-algebras. Its subalgebra
\[
G(A)=G_0\otimes A_0+G_1\otimes A_1
\]
is called the Grassmann envelope of the superalgebra~$A$.

Let $\Omega$ be a variety of algebras over~$F$. A superalgebra $A=A_0\oplus A_1$ is called an
$\Omega$-superalgebra if its Grassmann envelope $G(A)$ is an algebra in~$\Omega$. In particular,
$J=J_0\oplus J_1$ is a Jordan superalgebra if its Grassmann envelope $G(J)$ is a Jordan algebra. In
what follows, for a homogeneous element $x$ of a superalgebra $J=J_0\oplus J_1$ we put $p(x)=i$ if
$x\in J_i$. We denote the even part $J_0$ of a Jordan superalgebra by $A$ and the odd part $J_1$
by~$M$.

A classification of the simple finite-dimensional Jordan superalgebras over an algebraically closed
field of characteristic zero was given in~\cite{10,11}. In~\cite{12,13}, all simple finite-dimensional
Jordan superalgebras over an algebraically closed field of arbitrary characteristic different from~$2$
were described.

We give several examples of Jordan superalgebras.

\textbf{1.1. The Kantor double \cite{10}.} Let $\Gamma=\Gamma_0\oplus\Gamma_1$ be an associative
supercommutative superalgebra with unity~$1$, and let $\{\,,\}\colon\Gamma\times\Gamma\to\Gamma$ be a
superskewsymmetric even bilinear map, which we shall call a bracket. Starting with the
superalgebra $\Gamma$ and the bracket $\{\,,\}$, one can construct a superalgebra $\Js$. Consider the
following direct sum of spaces: $\Js=\Gamma\oplus\Gamma x$, where $\Gamma x$ is an isomorphic copy of
the space~$\Gamma$, $p(ax)=p(a)+1$. Let $a$ and $b$ be homogeneous elements in~$\Gamma$. Then
multiplication $\cdot$ on $\Js$ is defined by the relations
$$a\cdot b=ab,\qquad 
a\cdot bx=(ab)x,\qquad 
ax\cdot b=(-1)^{p(b)}(ab)x,\qquad
ax\cdot bx=(-1)^{p(b)}\br ab.$$
We set $A=\Gamma_0\oplus\Gamma_1x$ and $M=\Gamma_1\oplus\Gamma_0x$. Then $\Js=A\oplus M$ is a
$\mathbb Z_2$-graded algebra.

A bracket $\{\,,\}$ is said to be Jordan if the superalgebra $\Js$ is a Jordan superalgebra. It is well
known~\cite{14} that, for homogeneous elements, the Jordan bracket satisfies the relations
\begin{gather}
\br a{bc}=\br abc+(-1)^{p(a)p(b)}b\br ac-\br a1bc,\label{e1}\\
\begin{aligned}
\br a{\br bc}={}&\br{\br ab}c+(-1)^{p(a)p(b)}\br b{\br ac}+\br a1\br bc\\
&+(-1)^{p(a)(p(b)+p(c))}\br b1\br ca+(-1)^{p(c)(p(a)+p(b))}\br c1\br ab.
\end{aligned}\label{e2}
\end{gather}
Since the superalgebra $\Js$ is Jordan, we see that $D\colon a\mapsto\br a1$ is a derivation of the
superalgebra~$\Gamma$.

If $D$ is the zero derivation, then $\{\,,\}$ is a Poisson bracket, i.e.,
\[
\br a{bc}=\br abc+(-1)^{p(a)p(b)}b\br ac,
\]
and $\Gamma$ is a Lie superalgebra with respect to the operation $\{\,,\}$. An arbitrary Poisson
bracket is a Jordan bracket (see~\cite{15}).

It is well known~\cite{14} that the Jordan superalgebra $J=\Gamma\oplus\Gamma x$ obtained with the
help of the Kantor doubling process is simple if and only if $\Gamma$ has no nonzero
proper graded ideals $B$ satisfying $\{\Gamma,B\}\subseteq B$.

\textbf{1.2. A superalgebra of vector type $J(\Gamma,D)$.} Let $\Gamma=\Gamma_0\oplus\Gamma_1$ be an
associative supercommutative superalgebra with a nonzero even derivation~$D$. We define a bracket
$\{\,,\}$ on $\Gamma$ by setting
\[
\br ab=D(a)b-aD(b).
\]
Then the bracket $\{\,,\}$ is Jordan. We denote by $\Js$ the resulting superalgebra $J(\Gamma,D)$.
Multiplication ``$\cdot$'' in $J(\Gamma,D)$ is defined by
$$a\cdot b=ab,\quad a\cdot bx=(ab)x,\quad ax\cdot b=(-1)^{p(b)}(ab)x,\quad
ax\cdot bx=(-1)^{p(b)}(D(a)b-aD(b)),$$
where $a$ and $b$ are homogeneous elements in $\Gamma$ and $ab$ is the product in~$\Gamma$. The
superalgebra $J(\Gamma,D)$ is called a superalgebra of vector type. If the superalgebra $J(\Gamma,D)$
is simple, then $\Gamma_1=0$ (see~\cite{14}).

\textbf{1.3. The Cheng--Kac superalgebra $CK(Z,d)$ \cite{16}.} Let $Z$ be an arbitrary unital
associative commutative algebra with a nonzero derivation $d\colon Z\to Z$. Consider two free
$Z$-modules of rank~$4$:
\[
A=Z+\sum_{i=1}^3w_iZ,\qquad M=xZ+\sum_{i=1}^3x_iZ.
\]
Multiplication on $A$ will be $Z$-linear, $w_iw_j=0$, $i\neq j$, $w_i^2=-1$. We define
\[
x_{i\times i}=0,\quad x_{1\times2}=-x_{2\times1}=x_3,\quad
x_{1\times3}=-x_{3\times1}=-x_2,\quad x_{2\times3}=-x_{3\times2}=x_1.
\]
Multiplication $A\times M\to M$ is defined by the relations
\textup{(}$A$ and $M$ supercommute: $am=ma$\textup{)}
\[
(xf)g=x(fg),\quad (x_if)g=x_i(fg),\quad (xf)(w_jg)=x_j(fd(g)),\quad (x_if)(w_jg)=x_{i\times j}(fg).
\]
Multiplication $M\times M\to A$ is given in accordance with the rules
{\small $$
(xf)(xg)=d(f)g-fd(g),\quad (xf)(x_jg)=-w_j(fg),\quad (x_if)(xg)=w_i(fg),\quad (x_if)(x_jg)=0.
$$}

We also need the definition of a certain superalgebra $B(m,n)$. Let $F$ be an
algebraically closed field of characteristic $p>2$. Let $B(m)=F[a_1,\dots,a_m\mid a_i^p=0]$ be the
algebra of truncated polynomials in $m$ even variables. Let $G(n)$ be the Grassmann superalgebra with
generators $1,\xi_1,\dots,\xi_n$. Then $B(m,n)=B(m)\otimes G(n)$ is an associative supercommutative
superalgebra.

The main result concerning classification of simple finite-dimensional unital Jordan superalgebras
over algebraically closed fields of characteristic $p>2$ was obtained in the paper~\cite{13} by
Martinez and Zelmanov.

\begin{theorem}\label{thm1}
Let $J=J_0+J_1$ be a finite-dimensional simple unital Jordan superalgebra over an algebraically closed
field of characteristic $p>2$, where $J_0$ is not a semisimple algebra. Then\textup:

\begin{enumerate} 
\item[${\rm (1)}$]   either there exist integers $m$ and $n$ and a Jordan bracket $\{\,,\}$ on $B(m,n)$ such
that $J=J(B(m,n),\{\,,\})$,
    \item[${\rm (2)}$]  or $J$ is isomorphic to the Cheng--Kac Jordan superalgebra $CK(B(m),d)$, which is defined
by the derivation $d\colon B(m)\to B(m)$.

\end{enumerate}\end{theorem}

As has been mentioned above, for a fixed element $\delta$ of the ground field, by a
$\delta$-derivation of a superalgebra $A$ we mean a linear map $\varphi\colon A\to A$ such that
\[
\varphi(xy)=\delta(\varphi(x)y+x\varphi(y))
\]
for all $x,y\in A$.

The centroid $\Gamma(A)$ of a superalgebra $A$ is the set of all linear maps $\chi\colon A\to A$ such
that
\[
\chi(xy)=\chi(x)y=x\chi(y)
\]
for all $x,y\in A$.

Note that the $1$-derivation is the usual derivation and the $0$-derivation is an arbitrary
endomorphism $\varphi$ of the algebra $A$ such that $\varphi(A^2)=0$. It is clear that any element of
the centroid of an algebra is a $\frac12$-derivation.

A nonzero $\delta$-derivation $\varphi$ is said to be nontrivial if $\delta\notin\{0,1\}$ and
$\varphi\notin\Gamma(A)$.

By a superspace we mean a $\mathbb Z_2$-graded space. A homogeneous element $\psi$ of the superspace
of endomorphisms $A\to A$ is called a superderivation if
\[
\psi(xy)=\psi(x)y+(-1)^{p(x)p(\psi)}x\psi(y).
\]

For a fixed element $\delta\in F$, we define the notion of a $\delta$-superderivation of the
superalgebra $A=A_0+A_1$. A homogeneous linear map $\varphi\colon A\to A$ is called a
$\delta$-superderivation if for homogeneous $x,y\in A$ we have
\[
\varphi(xy)=\delta(\varphi(x)y+(-1)^{p(x)p(\varphi)}x\varphi(y)).
\]
Consider a Lie superalgebra $A=A_0+A_1$ and fix an element $x\in A_i$. Then $R_x\colon y\mapsto xy$
is a superderivation of the superalgebra $A$ and its parity $p(R_x)$
is equal to~$i$.

The supercentroid $\Gamma_s(A)$ of a superalgebra $A$ is the set of all homogeneous linear maps
$\chi\colon A\to A$ such that for arbitrary homogeneous elements $a$ and $b$, we have
\[
\chi(ab)=\chi(a)b=(-1)^{p(a)p(\chi)}a\chi(b).
\]
Note that a $1$-superderivation is a usual superderivation; a $0$-superderivation is an arbitrary
homogeneous endomorphism $\varphi$ of the superalgebra $A$ such that $\varphi(A^2)=0$.

A nonzero $\delta$-superderivation $\varphi$ is said to be nontrivial if $\delta\notin\{0,1\}$ and
$\varphi\notin\Gamma_s(A)$.

In accordance with~\cite[Theorem~2.1]{4} (which is easily generalized to the case of
$\delta$-superderivations), for a unital superalgebra $A$, a map $\varphi$ can be a nontrivial
$\delta$-derivation or $\delta$-superderivation only for $\delta=\frac12$. It is easily seen that, in
this case, $\varphi(x)=\varphi(1)x$ for arbitrary $x\in A$.

\section{$\delta$-derivations and $\delta$-superderivations of simple superalgebras of Jordan brackets}

In the present section, we consider $\delta$-derivations and $\delta$-superderivations of a simple
unital Jordan superalgebra $J=\Js$. We assume that the characteristic of the field $F$ is different
from~$2$.

\begin{lemma}\label{l2}
Let $J=\Js$ be a simple unital Jordan superalgebra. Then $\Gamma=\Gamma\{\Gamma,\Gamma\}$. In
particular, if $z\in\Gamma_0\cup\Gamma_1\setminus\{0\}$, then $z\{\Gamma,\Gamma\}\neq0$.
\end{lemma}

\begin{proof}
Consider $I=\Gamma\{\Gamma,\Gamma\}$. It is clear that $I$ is an ideal in $\Gamma$ ($I\trianglelefteq\Gamma$).
By~\eqref{e1},
\[
\{\Gamma,I\}=\{\Gamma,\Gamma\{\Gamma,\Gamma\}\}\subseteq
\{\Gamma,\Gamma\}\{\Gamma,\Gamma\}
+\Gamma\{\Gamma,\{\Gamma,\Gamma\}\}+\{\Gamma,1\}\Gamma\{\Gamma,\Gamma\}\subseteq\Gamma\{\Gamma,\Gamma\}=I.
\]
By~\cite{14}, the Jordan superalgebra $\Js$ is simple if $\Gamma$ contains no nonzero proper
ideals $I$ satisfying the condition $\{\Gamma,I\}\subseteq I$. If $\{\Gamma,\Gamma\}=0$, then
$\Gamma x\trianglelefteq J$. Consequently, $\Gamma x=J$, which is
impossible, since $1\notin\Gamma x$. Hence, $\{\Gamma,\Gamma\}\neq0$, and since $\Gamma$ is unital,
we have $\Gamma=\Gamma\{\Gamma,\Gamma\}$. Now if $z\{\Gamma,\Gamma\}=0$, then
$z\Gamma=z\Gamma\{\Gamma,\Gamma\}=0$. Since $\Gamma$ is unital, we obtain $z=0$. The lemma is proved.
\end{proof}

\begin{lemma}\label{l3}
Let $J=\Js$ be a simple unital Jordan superalgebra, and let $\alpha\in J$. The map
$\varphi(z)=\alpha z$ is a $\frac12$-derivation if and only if $\alpha\in\Gamma_0$ and
$\br\alpha b=D(\alpha)b-\alpha D(b)$ for any $b\in\Gamma$.
\end{lemma}

\begin{proof}
Let $\alpha=\alpha_0+\beta x+\gamma+\mu x$, where $\alpha_0,\mu\in\Gamma_0$ and
$\beta,\gamma\in\Gamma_1$. Clearly, the maps $\varphi_1(z)=(\gamma+\mu x)z$ and
$\varphi_2(z)=(\alpha_0+\beta x)z$ are also $\frac12$-derivations of the superalgebra~$J$. Moreover,
$\varphi_1(1)=(\gamma+\mu x)$ and $\varphi_2(1)=(\alpha_0+\beta x)$. For this reason, for arbitrary
$z,w\in J$ we have
\begin{equation}
2\varphi_i(1)(zw)=(\varphi_i(1)z)w+z(\varphi_i(1)w).\label{e3}
\end{equation}
Setting $i=1$, $z=x$, and $w=1$ in~\eqref{e3}, we get
\[
2\gamma x=\gamma x+x\gamma=0,
\]
i.e., $\gamma=0$.

We prove that $\beta=\mu=0$. For this, we show that $\beta\{\Gamma,\Gamma\}=0$ and
$\mu\{\Gamma,\Gamma\}=0$.

Now in~\eqref{e3} we set $i=2$, $z=a$, and $w=bx$, obtaining
\[
2\br\beta{ab}=\br{\beta a}b+(-1)^{p(a)}a\br\beta b.
\]
By~\eqref{e1}, we have
\begin{equation}
2\br\beta{ab}=-(-1)^{p(b)+p(a)p(b)}2\br b\beta a-(-1)^{p(a)p(b)}\beta\br ba
+(-1)^{p(b)+p(b)p(a)}D(b)\beta a.\label{e4}
\end{equation}
Substituting $i=2$, $z=ax$, and $w=b$ in~\eqref{e3}, we get
\[
2\br\beta{ab}=\br\beta ab-(-1)^{p(a)}\br a{\beta b}.
\]
By~\eqref{e1}, we have
\begin{equation}
2\br\beta{ab}=2\br\beta ab-\beta\br ab-\beta D(a)b.\label{e5}
\end{equation}
Substituting $i=2$, $z=ax$, and $w=bx$ in~\eqref{e3}, we get
\begin{equation}
2\beta\br ab=\br\beta ab-(-1)^{p(a)}a\br\beta b.\label{e6}
\end{equation}
Comparison of~\eqref{e4} and~\eqref{e5} yields
\begin{equation}
\begin{aligned}
2\beta\br ab&=2\br\beta ab-\beta D(a)b-(-1)^{p(b)+p(b)p(a)}D(b)\beta a+(-1)^{p(b)+p(b)p(a)}2\br b\beta a\\
&=2\br\beta ab-\beta D(a)b-\beta aD(b)-(-1)^{p(a)}2a\br\beta b\\
&=2\br\beta ab-\beta D(ab)-(-1)^{p(a)}2a\br\beta b.
\end{aligned}\label{e7}
\end{equation}
This relation and~\eqref{e6} imply that
\[
2\beta\br ab=\beta D(ab).
\]
Putting $b=1$, we get $\beta D(a)=2\beta D(a)$. Hence, $\beta D(a)=0$ and
$\beta\{\Gamma,\Gamma\}=0$. Thus, $\beta=0$ by Lemma~\ref{l2}.

Substituting $i=1$, $z=a$, $w=bx$ in~\eqref{e3}, we obtain
\begin{equation}
2\br\mu{ab}=-\br{\mu a}b+(-1)^{p(a)}a\br\mu b.\label{e8}
\end{equation}
Choosing $i=1$, $z=ax$, $w=b$ in~\eqref{e3}, we get
\begin{equation}
2\br\mu{ab}=\br\mu ab+(-1)^{p(a)}\br a{\mu b}.\label{e9}
\end{equation}
For $i=1$, $z=ax$, $w=bx$ in~\eqref{e3}, we have
\begin{equation}
2\mu\br ab=\br\mu ab+(-1)^{p(a)}a\br\mu b.\label{e10}
\end{equation}
For $a=b=1$, relation~\eqref{e10} implies that $D(\mu)=0$. Substituting $b=1$ in~\eqref{e8}, we
obtain
\[
2\br\mu a=\br{\mu a}1=D(\mu a)=\mu D(a).
\]
Taking $a=1$ in~\eqref{e9}, we have
\[
2\br\mu b=\br1{\mu b}=-D(\mu b)=-\mu D(b).
\]
Comparing the expressions obtained, we get $\{\mu,\Gamma\}=0$. Therefore, $\mu\br ab=0$
by~\eqref{e10}. Using Lemma~\ref{l2}, we see that $\mu=0$.

Thus, we have proved that $\varphi(z)=\alpha z$, where $\alpha\in\Gamma_0$.

Putting $i=2$, $z=ax$, and $w=bx$ in~\eqref{e3}, we obtain
\begin{equation}
2\alpha\br ab=\br{\alpha a}b+\br a{\alpha b}.\label{e11}
\end{equation}
Note that identity~\eqref{e1} yields the relation
\begin{multline*}
\br{\alpha a}b+\br a{\alpha b}=\\-(-1)^{p(b)p(a)}\bigl(\br b\alpha a+\alpha\br ba-D(b)\alpha a\bigr)
+\br a\alpha b+\alpha\br ab-D(a)\alpha b.
\end{multline*}
Hence, \eqref{e11} shows that
\[
a\br b\alpha-\br a\alpha b=(aD(b)-D(a)b)\alpha.
\]
Consequently, for $b=1$,
\[
\br\alpha a=D(\alpha)a-\alpha D(a).
\]
It is easy to verify that, for any $\alpha\in\Gamma_0$ such that $\br\alpha a=D(\alpha)a-\alpha D(a)$,
$\varphi(z)=\alpha z$ is a $\frac12$-derivation of the superalgebra~$J$. The lemma is proved.
\end{proof}

Thus, a $\delta$-derivation of a simple unital superalgebra $J=\Js$ is an even
$\delta$-superderivation of the superalgebra $J=\Js$. For this reason, in the sequel we deal only with
$\delta$-superderivations.

\begin{remark}\label{r4}
Let $J=\Js$ be a simple unital Jordan superalgebra over a field of characteristic $p\neq2,3$.
The map $\varphi(z)=\alpha z$ is an odd $\frac12$-superderivation if and only if $\alpha\in\Gamma_1$
and
\[
\br\alpha a=D(\alpha)a-\alpha D(a)
\]
for any $a\in\Gamma$. If $p=3$, the odd $\frac12$-superderivations of the form $\varphi(z)=\alpha z$
are the maps with $\alpha=\gamma+\mu x$, where $\gamma$ is as above and $\mu\in\Gamma_0$ satisfies a
system of conditions which implies $\mu\in\Phi$ \textup{(}see below\textup{)} and
$\mu\br ab=\mu(D(a)b-aD(b))$ for all $a,b\in\Gamma$; if $J=J(\Gamma,D)$ is of vector type, then every
$\mu\in\Gamma$ is admissible \textup{(}Remark~A\textup{)}.
\end{remark}

\begin{proof}
This is proved by straightforward calculations similar to those in the proof of Lemma~\ref{l3}
\textup{(}the component $\mu x$ of $\alpha$ contributes the relations $3\mu D(a)=0$; in
characteristic $\neq3$ they give $\mu=0$ as in Lemma~\ref{l3}; see Remark~A\textup{)}.
\end{proof}

\begin{corollary}\label{c5}
If $J$ is a simple unital superalgebra of vector type, then the map $\varphi(z)=\alpha z$ is a
$\frac12$-superderivation if and only if $\alpha\in\Gamma_0$ \textup{(}if $p\neq3$\textup{)}, resp.
$\alpha\in\Gamma_0\cup\Gamma_0x$ \textup{(}if $p=3$\textup{)}. If $J$ is a superalgebra of the Poisson
bracket, then the map $\varphi(z)=\alpha z$ is a $\frac12$-superderivation if and only if
$\alpha\in\Gamma_0\cup\Gamma_1$ and $\{\alpha,\Gamma\}=0$.
\end{corollary}

Let $J=J(\Gamma,D)$ be a superalgebra of vector type. The map $\varphi(z)=\alpha z$ with
$\alpha\in\Gamma$ is the trivial $\frac12$-superderivation if $\varphi\in\Gamma_s(J)$, i.e., if
\[
\alpha((bx)(cx))=(-1)^{p(\alpha)p(bx)}(bx)(\alpha(cx)),
\]
which is equivalent to $D(\alpha)bc=0$. Consequently, $\varphi$ is the nontrivial
$\frac12$-superderivation if $D(\alpha)\neq0$.

Let $J=\Js$ be a superalgebra of the Poisson bracket, and let $\varphi(z)=\alpha z$ be a
$\frac12$-superderivation of~$J$. By Remark~\ref{r4}, we have
\begin{align*}
(ax)(\alpha(bx))&=(-1)^{p(b)+p(\alpha)}\br a{\alpha b}
=(-1)^{p(b)+p(\alpha)}\bigl(\br a\alpha b+(-1)^{p(\alpha)p(a)}\alpha\br ab\bigr)\\
&=(-1)^{p(b)+p(\alpha)+p(a)p(\alpha)}\alpha\br ab=(-1)^{p(\alpha)(p(a)+1)}\alpha((ax)(bx)).
\end{align*}
This easily implies that $\varphi$ is the trivial $\frac12$-superderivation.

We set
\[
\Phi=\{\alpha\in\Gamma_0\cup\Gamma_1\mid \br\alpha a=D(\alpha)a-\alpha D(a),\ a\in\Gamma\}.
\]

\begin{lemma}\label{l6}
Suppose $J=\Js$ is a simple Jordan superalgebra, $\alpha\in\Phi\setminus\{0\}$, and $D(\alpha)=0$.
Then $\alpha$ is invertible in $\Gamma$ and $\alpha\in\Gamma_0$. In particular, if $\Js$ is a
superalgebra of the Poisson bracket, then $\alpha$ is invertible in $\Gamma$ and $\alpha\in\Gamma_0$.
\end{lemma}

\begin{proof}
Consider $I=\alpha\Gamma$. It is clear that $I\trianglelefteq\Gamma$. By the definition of $\Phi$, we
have
\[
\{I,\Gamma\}=\{\alpha\Gamma,\Gamma\}\subseteq\{\Gamma,\alpha\}\Gamma+\alpha\{\Gamma,\Gamma\}
+D(\Gamma)\alpha\Gamma\subseteq\alpha\Gamma=I.
\]
Using~\cite{14}, we conclude that $I=\Gamma$. Since $\Gamma$ is unital, it follows that $\alpha$ is
invertible \textup{(}and hence even, since odd elements of $\Gamma$ have zero square\textup{)}.
The lemma is proved.
\end{proof}

\begin{lemma}\label{l7}
Let $J=\Js$ be a Jordan superalgebra. Then $\Phi$ is closed with respect to the derivation $D$, i.e.,
$D(\Phi)\subseteq\Phi$. In particular, $D^k(\Phi)\subseteq\Phi$, $k>0$.
\end{lemma}

\begin{proof}
Using~\eqref{e2}, we obtain
\begin{align*}
\br{D(b)}c+\br b{D(c)}&=-(-1)^{p(b)p(c)}\br c{\br b1}+\br b{\br c1}\\
&=-(-1)^{p(b)p(c)}\br{\br cb}1-(-1)^{p(b)p(c)}(-1)^{p(c)p(b)}\br b{\br c1}\\
&\quad-(-1)^{p(b)p(c)}\br c1\br b1-(-1)^{p(b)p(c)}(-1)^{p(c)p(b)}\br b1\br1c\\
&\quad+\br{\br bc}1+(-1)^{p(b)p(c)}\br c{\br b1}\\
&\quad+\br b1\br c1+(-1)^{p(b)p(c)}\br c1\br1b\\
&=2D(\br bc)-\br{D(b)}c-\br b{D(c)}.
\end{align*}
Consequently,
\begin{equation}
D(\br bc)=\br{D(b)}c+\br b{D(c)}.\label{e12}
\end{equation}
Therefore,
\begin{align*}
\br{D(\alpha)}a&=D(\br\alpha a)-\br\alpha{D(a)}\\
&=D(D(\alpha))a+D(\alpha)D(a)-D(\alpha)D(a)-\alpha D(D(a))-D(\alpha)D(a)+\alpha D(D(a))\\
&=D(D(\alpha))a-D(\alpha)D(a).
\end{align*}
This completes the proof of the lemma.
\end{proof}

\begin{lemma}\label{l8}
Let $J=\Js$ be a Jordan superalgebra. Then for arbitrary $b,c\in\Gamma$ and
$\alpha\in\Phi$ we have
\[
D^k(\alpha)\br bc=D^k(\alpha)(D(b)c-bD(c)).
\]
\end{lemma}

\begin{proof}
From the definition of $\Phi$ and~\eqref{e12}, it follows that
\begin{equation}
\br\alpha{\br bc}=D(\alpha)\br bc-\alpha D(\br bc)=D(\alpha)\br bc-\alpha\br{D(b)}c-\alpha\br b{D(c)}.
\label{e13}
\end{equation}
Now we use relations~\eqref{e1}, \eqref{e2} and Lemmas~\ref{l3} and~\ref{l7} to write
\begin{longtable}{lcl}
$\br\alpha{\br bc}$&$=$&
$\br{\br\alpha b}c+(-1)^{p(\alpha)p(b)}\br b{\br\alpha c}+D(\alpha)\br bc$\\
&&$+(-1)^{p(\alpha)(p(b)+p(c))}D(b)\br c\alpha+(-1)^{p(c)(p(\alpha)+p(b))}D(c)\br\alpha b$\\
&$=$&$\br{D(\alpha)b}c-\br{\alpha D(b)}c+(-1)^{p(\alpha)p(b)}\br b{D(\alpha)c}
-(-1)^{p(\alpha)p(b)}\br b{\alpha D(c)}$\\
&&$+D(\alpha)\br bc-(-1)^{p(\alpha)p(b)}D(b)(D(\alpha)c-\alpha D(c))$\\
&&$+(-1)^{p(c)(p(\alpha)+p(b))}D(c)(D(\alpha)b-\alpha D(b))$\\
&$=$&$-(-1)^{p(c)(p(\alpha)+p(b))}\br c{D(\alpha)b}+(-1)^{p(c)(p(b)+p(\alpha))}\br c{\alpha D(b)}$\\
&&$+(-1)^{p(\alpha)p(b)}\br b{D(\alpha)c}-(-1)^{p(\alpha)p(b)}\br b{\alpha D(c)}+D(\alpha)\br bc$\\
&&$-(-1)^{p(\alpha)p(b)}D(b)D(\alpha)c+(-1)^{p(\alpha)p(b)}D(b)\alpha D(c)$\\
&&$+(-1)^{p(c)(p(\alpha)+p(b))}D(c)D(\alpha)b-(-1)^{p(c)(p(\alpha)+p(b))}D(c)\alpha D(b)$\\
&$=$&
$-(-1)^{p(c)(p(b)+p(\alpha))}\br c{D(\alpha)}b-(-1)^{p(b)p(c)}D(\alpha)\br cb$\\
&&$+(-1)^{p(c)(p(\alpha)+p(b))}D(c)D(\alpha)b
+(-1)^{p(c)(p(b)+p(\alpha))}\br c\alpha D(b)$\\
&&$+(-1)^{p(b)p(c)}\alpha\br c{D(b)}
+(-1)^{p(c)(p(\alpha)+p(b))}D(c)\alpha D(b)$\\
&&$+(-1)^{p(\alpha)p(b)}\br b{D(\alpha)}c+D(\alpha)\br bc-(-1)^{p(\alpha)p(b)}D(b)D(\alpha)c$\\
&&$-(-1)^{p(\alpha)p(b)}\br b\alpha D(c)-\alpha\br b{D(c)}+(-1)^{p(\alpha)p(b)}D(b)\alpha D(c)$\\
&&$+D(\alpha)\br bc
-(-1)^{p(\alpha)p(b)}D(b)D(\alpha)c+(-1)^{p(\alpha)p(b)}D(b)\alpha D(c)$\\ 
&&$+(-1)^{p(c)(p(\alpha)+p(b))}D(c)D(\alpha)b
-(-1)^{p(c)(p(\alpha)+p(b))}D(c)\alpha D(b)$\\
&$=$&$(-1)^{p(b)p(c)}D(D(\alpha))cb-(-1)^{p(c)p(b)}D(\alpha)D(c)b+D(\alpha)\br bc$\\
&&$+(-1)^{p(c)(p(\alpha)+p(b))}D(c)D(\alpha)b-(-1)^{p(b)p(c)}D(\alpha)cD(b)$\\
&&$+(-1)^{p(b)p(c)}\alpha D(c)D(b)-\alpha\br{D(b)}c+(-1)^{p(c)(p(\alpha)+p(b))}D(c)\alpha D(b)$\\
&&$-D(D(\alpha))bc+D(\alpha)D(b)c+D(\alpha)\br bc -(-1)^{p(\alpha)p(b)}D(b)D(\alpha)c$\\
&&$+D(\alpha)bD(c) -\alpha D(b)D(c)-\alpha\br b{D(c)} +(-1)^{p(\alpha)p(b)}D(b)\alpha D(c)$\\
&&$+D(\alpha)\br bc-(-1)^{p(\alpha)p(b)}D(b)D(\alpha)c +(-1)^{p(\alpha)p(b)}D(b)\alpha D(c)$\\
&&$+(-1)^{p(c)(p(\alpha)+p(b))}D(c)D(\alpha)b
-(-1)^{p(c)(p(\alpha)+p(b))}D(c)\alpha D(b)$\\
&$=$&$3D(\alpha)\br bc-\alpha\br{D(b)}c-\alpha\br b{D(c)}-2D(\alpha)(D(b)c-bD(c)).$
\end{longtable}
Therefore,
\begin{equation}
3D(\alpha)\br bc-\alpha\br{D(b)}c-\alpha\br b{D(c)}-2D(\alpha)(D(b)c-bD(c))=\br\alpha{\br bc}.
\label{e14}
\end{equation}
Comparing~\eqref{e13} and~\eqref{e14}, we obtain
\[
D(\alpha)\br bc=D(\alpha)(D(b)c-bD(c)).
\]
Lemma~\ref{l7} allows us to generalize this relation as follows:
\[
D^k(\alpha)\br bc=D^k(\alpha)(D(b)c-bD(c)).
\]
The lemma is proved.
\end{proof}

\begin{lemma}\label{l9}
Let $J=\Js$ be a simple unital Jordan superalgebra, and let $\alpha\in\Phi$. If $D(\alpha)\neq0$, then
$J$ is a superalgebra of vector type. In particular, if $J$ is not a superalgebra of vector type, then
$D(\alpha)=0$ and for $\alpha\neq0$ $\alpha$ is invertible in~$\Gamma$.
\end{lemma}

\begin{proof}
We set $I=\Gamma D(\alpha)+\Gamma D^2(\alpha)+\cdots$. Note that $I\trianglelefteq\Gamma$.
Using~\eqref{e1}, the definition of $\Phi$, and Lemma~\ref{l7}, it is easy to deduce that
\begin{align*}
\{\Gamma,\Gamma D^k(\alpha)\}&\subseteq\{\Gamma,D^k(\alpha)\}\Gamma+D^k(\alpha)\{\Gamma,\Gamma\}
+D(\Gamma)D^k(\alpha)\Gamma\\
&\subseteq D(\Gamma)D^k(\alpha)\Gamma+\Gamma D^{k+1}(\alpha)\Gamma+D^k(\alpha)\{\Gamma,\Gamma\}
+D(\Gamma)D^k(\alpha)\Gamma\\
&\subseteq\Gamma D^k(\alpha)+\Gamma D^{k+1}(\alpha),
\end{align*}
whence $\{\Gamma,I\}\subseteq I$. By~\cite{14}, $\Gamma$ contains no nonzero proper ideals $I$
such that $\{\Gamma,I\}\subseteq I$. If $D(\alpha)\neq0$, then $I=\Gamma$. Therefore,
$1=\gamma_1D(\alpha)+\dots+\gamma_lD^l(\alpha)$. Consequently, by Lemma~\ref{l8}, for arbitrary
$b,c\in\Gamma$ we get
\begin{align*}
\br bc&=(\gamma_1D(\alpha)+\dots+\gamma_lD^l(\alpha))\br bc\\
&=\gamma_1D(\alpha)(D(b)c-bD(c))+\dots+\gamma_lD^l(\alpha)(D(b)c-bD(c))=D(b)c-bD(c).
\end{align*}
Thus, $J$ is a superalgebra of vector type.

If $J$ is not a superalgebra of vector type, then the said above implies $D(\alpha)=0$. By
Lemma~\ref{l6}, a nonzero $\alpha$ is an invertible element in~$\Gamma$. The lemma is proved.
\end{proof}

The above results are generalized in the following theorem \textup{(}for a proof see
Remark~B\textup{)}.

\begin{theorem}\label{t10}
Let $J=\Js$ be a simple unital superalgebra of a Jordan bracket over a field of characteristic
different from~$2$. Then either $J$ has no nontrivial $\delta$-derivations and
$\delta$-superderivations, or $J$ is a superalgebra of vector type. If $J$ is a superalgebra of vector
type, then $\Gamma_1=0$ and for $p\neq3$ the superalgebra $J$ has no nontrivial odd
$\delta$-superderivations; for $p=3$ the odd $\delta$-superderivations are zero for
$\delta\neq\frac12$, and the odd $\frac12$-superderivations are exactly the maps
$R_{\mu x}\colon z\mapsto(\mu x)z$, $\mu\in\Gamma$, which are nontrivial for $\mu\neq0$. For
$\delta\neq\frac12$, the superalgebra $J$ has no nontrivial $\delta$-derivations. The space of
$\frac12$-derivations coincides with $R^*(J)=\{R_z\mid z\in\Gamma_0\}$, and if $D(z)\neq0$, then the
map $R_z$ is a nontrivial $\frac12$-derivation.
\end{theorem}

\begin{RemA}
Let $p=3$ and let $J=\Js$ be a simple unital Jordan superalgebra.

(a) \emph{Necessary conditions.} Let $\alpha=\gamma+\mu x$, $\gamma\in\Gamma_1$, $\mu\in\Gamma_0$,
and let $\varphi(z)=\alpha z$ be an odd $\frac12$-superderivation. The multiplication of $J$ respects
the decomposition $J=\Gamma\oplus\Gamma x$ modulo $2$, so $z\mapsto\gamma z$ and $z\mapsto(\mu x)z$
are odd $\frac12$-superderivations separately. For $z\mapsto(\mu x)z$ the identity
$\varphi(zw)=\frac12(\varphi(z)w+(-1)^{p(z)}z\varphi(w))$ at $(z,w)=(a,x)$, $(ax,x)$, $(ax,bx)$ gives
$$\br\mu a=D(\mu)a+\tfrac12\mu D(a),\qquad \br\mu a=aD(\mu)+2\mu D(a),\\
\mu\br ab=\tfrac12\bigl(\br\mu ab-a\br\mu b\bigr).$$

The first two relations give $3\mu D(a)=0$; for $p\neq3$ this forces $\mu D(\Gamma)=0$ and, as in
the proof of Lemma~\ref{l3}, $\mu=0$. For $p=3$ they coincide and mean $\mu\in\Phi$; then the third one
becomes $\mu\br ab=-\frac12\mu(D(a)b-aD(b))=\mu(D(a)b-aD(b))$.

(b) \emph{$J$ not of vector type.} If $\mu\neq0$, then by Lemma~\ref{l9} $D(\mu)=0$ and $\mu$ is
invertible, and $\mu\br ab=\mu(D(a)b-aD(b))$ shows that $\br ab=D(a)b-aD(b)$, i.e., $J$ is of vector
type, a contradiction. Hence $\mu=0$, and the odd $\frac12$-superderivations are as in Remark~\ref{r4}.

(c) \emph{$J=J(\Gamma,D)$ of vector type.} Here $\Gamma=\Gamma_0$. A direct computation shows
that for every $\mu\in\Gamma$ the map $R_{\mu x}\colon a+bx\mapsto\br\mu b+(\mu a)x$ is an odd
$\frac12$-superderivation: the defining identity at $(a,c)$, $(a,dx)$, $(bx,c)$ reduces to 
$$3\mu ac=0,\quad
3(D(\mu)a d-\mu aD(d))=0,\quad 
3(D(\mu)bc-\mu cD(b))=0,$$
and at $(bx,dx)$ it holds identically. If
$R_{\mu x}\in\Gamma_s(J)$, then $R_{\mu x}((1x)(ax))=R_{\mu x}(1x)\cdot ax$ gives $-\mu D(a)=D(\mu)a$
for all $a$, whence $D(\mu)=0$ and $\mu D(\Gamma)=0$; then $\{\Gamma,\mu\Gamma\}\subseteq\mu\Gamma$, so
$\mu\Gamma=\Gamma$ for $\mu\neq0$, $\mu$ is invertible and $D=0$, a contradiction. Thus $R_{\mu x}$ is
nontrivial for every $\mu\neq0$ \textup{(}already $\mu=1$\textup{)}.

All statements of~(c), the relation $3\mu D(a)=0$ in every characteristic and the nontriviality of
$R_{\mu x}$ are verified in Lean.
\end{RemA}

\begin{RemB}
\emph{Proof of Theorem~\ref{t10}.} Since $J$ is unital, a nontrivial $\delta$-(super)derivation
has $\delta=\frac12$ and is of the form $\varphi(z)=\alpha z$ with $\alpha=\varphi(1)$ homogeneous.
If $\varphi$ is even, then by Lemma~\ref{l3} $\alpha\in\Phi\cap\Gamma_0$; if $\varphi$ is odd, then by
Remark~\ref{r4} and Remark~A $\alpha=\gamma+\mu x$ with $\gamma\in\Phi\cap\Gamma_1$. If $J$ is not of
vector type, then by Lemma~\ref{l9} every nonzero element of $\Phi$ is invertible and has $D=0$; thus
$\gamma=0$ \textup{(}odd elements are not invertible\textup{)}, $\mu=0$ by Remark~A(b), and for even
$\varphi$ the element $\alpha$ satisfies $D(\alpha)=0$ and $\br\alpha a=-\alpha D(a)$; then
$z\mapsto\alpha z$ lies in the supercentroid \textup{(}this follows from~\eqref{e1} as in the
computation for the Poisson bracket above\textup{)}, i.e., $\varphi$ is trivial. If $J$ is of vector type, then
$\Gamma_1=0$~\cite{14}, the even case is Corollary~\ref{c5} and the paragraph following it, and the
odd case is Remark~\ref{r4} \textup{(}$p\neq3$\textup{)} and Remark~A(c) \textup{(}$p=3$\textup{)}.
\end{RemB}

Let $A=A_0\oplus A_1$ be an associative superalgebra. On the vector space $A$, we define a
supersymmetric product $\circ$ by the rule
\[
a\circ b=\tfrac12\bigl(ab+(-1)^{p(a)p(b)}ba\bigr).
\]
We denote by $A^{(+)}$ the algebra obtained. A Jordan superalgebra $B$ is said to be special if it is
isomorphically embedded in the superalgebra $A^{(+)}$ for a suitable associative superalgebra~$A$.

Using Theorem~\ref{t10} and the well-known fact that the unital Jordan superalgebras of vector type are
special (see~\cite{14}), we obtain the following statement.

\begin{corollary}\label{c11}
If a simple unital superalgebra of a Jordan bracket $J=\Js$ has a nontrivial $\delta$-derivation, then
$J$ is special.
\end{corollary}

\section{The $\delta$-derivations and $\delta$-superderivations of simple unital finite-dimensional Jordan superalgebras}

Now we proceed to a description of the $\delta$-derivations and $\delta$-superderivations of simple
unital finite-dimensional Jordan superalgebras over an algebraically closed field of characteristic
$p\neq2$. We recall that the $\delta$-derivations and $\delta$-superderivations of simple unital
finite-dimensional Jordan superalgebras over an algebraically closed field of characteristic zero were
described in~\cite{4,6}.

Recall that an algebra $A$ is said to be alternative if $(x,x,y)=0$ and $(x,y,y)=0$, where $(x,y,z)=(xy)z-x(yz)$ is the associator of elements $x,y,z\in A$. The algebra $\mathbb O$
of octonions or Cayley numbers (see~\cite{17}) is a classical example of an alternative nonassociative
algebra.

Now we give examples of simple nontrivial nonassociative alternative superalgebras of characteristic~$3$.
Below $B$ denotes an alternative superalgebra over a field $F$, and $C$ and $M$ are the even and the odd
part of $B$, respectively.

\textbf{3.1. The superalgebra $B(1,2)$.} Let $F$ be a field of characteristic~$3$, and let
$B(1,2)=C+M$ be a supercommutative superalgebra over $F$ in which $C=F\cdot1$ and $M=F\cdot x+F\cdot y$,
where $1$ is the unity of $B$ and $xy=-yx=1$. Note that the superalgebra $B(1,2)$ is precisely the
simple Jordan superalgebra of the supersymmetric bilinear form $f(s,r)=sr$ on the odd vector space~$M$.

\textbf{3.2. The superalgebra $B(4,2)$.} Let $F$ be a field of characteristic~$3$, let $C=M_2(F)$ be
the algebra of $(2\times2)$-matrices over $F$, and let $M=F\cdot m_1+F\cdot m_2$ be a two-dimensional
irreducible Cayley bimodule over $C$, i.e., $C$ acts on $M$ in the following way:
\[
e_{ij}\cdot m_k=\delta_{ik}m_j,\quad i,j,k\in\{1,2\},\qquad m\cdot a=\bar a\cdot m,
\]
where $a\in C$, $m\in M$, and $a\mapsto\bar a$ is the symplectic involution in $C=M_2(F)$. Odd
multiplication on $M$ is defined by the relations
\[
m_1^2=-e_{21},\qquad m_2^2=e_{12},\qquad m_1m_2=e_{11},\qquad m_2m_1=-e_{22}.
\]
It is known (see~\cite{18}) that $B(1,2)$ and $B(4,2)$ are simple alternative superalgebras with the
superinvolutions
\[
(a+m)^*=a-m\ \text{ for } B(1,2),\qquad (a+m)^*=\bar a-m\ \text{ for }
B(4,2).
\]
An even linear transformation $*$ of a superalgebra $A$ is called a superinvolution if
\[
(a^*)^*=a,\qquad (ab)^*=(-1)^{p(a)p(b)}b^*a^*,\qquad a,b\in A_0\cup A_1.
\]

In~\cite{18} it was shown that the superalgebras $B(1,2)$ and $B(4,2)$ give rise to the simple Jordan
superalgebras $H_3(B(1,2))$ and $H_3(B(4,2))$.

\begin{lemma}\label{l12}
The superalgebras $H_3(B(1,2))$ and $H_3(B(4,2))$ have no nontrivial
$\delta$-derivations and $\delta$-superderivations.
\end{lemma}

\begin{proof}
We denote by $e_{ij}$ the matrix units of the algebras $B(1,2)_3$ and $B(4,2)_3$. Let
$\varphi$ be a nontrivial $\delta$-derivation or $\delta$-superderivation. Clearly, $\delta=\frac12$.
Let
\[
\varphi(e_{ii})=\sum_{j=1}^3\alpha^i_je_{jj}+\sum_{k,l,\,k\neq l}x^i_{kl}e_{kl},
\]
where $x^i_{kl}=\bar x^i_{lk}$\textup{(}we write $x\,e_{kl}+\bar x\,e_{lk}$
for the element $x[kl]$ of $H_3(B)$\textup{)}; then
\[
2\varphi(e_{ii})=2e_{ii}\circ\varphi(e_{ii})=2\alpha^i_ie_{ii}+\sum_{k\neq i}(x^i_{ik}e_{ik}+x^i_{ki}e_{ki}).
\]
Thus, $\varphi(e_{ii})=\alpha^i_ie_{ii}$. If $\beta\in F$ and
$
\varphi(\beta(e_{21}+e_{12}))=\begin{pmatrix}\gamma_1&x_{12}&x_{13}\\ x_{12}&\gamma_2&x_{23}\\
x_{13}&x_{23}&\gamma_3\end{pmatrix},$
then
\[
\begin{pmatrix}0&x_{12}+\frac12\beta\alpha_{22}&0\\ x_{12}+\frac12\beta\alpha_{22}&2\gamma_2&x_{23}\\
0&x_{23}&0\end{pmatrix}=\varphi(\beta(e_{12}+e_{21}))\circ e_{22}+\beta(e_{12}+e_{21})\circ\varphi(e_{22})
\]
\[
=2\varphi\bigl(\beta(e_{12}+e_{21})\circ e_{22}\bigr)=2\varphi\bigl(\beta(e_{12}+e_{21})\circ e_{11}\bigr)
=\begin{pmatrix}2\gamma_1&x_{12}+\frac12\beta\alpha_{11}&x_{13}\\ x_{12}+\frac12\beta\alpha_{11}&0&0\\
x_{13}&0&0\end{pmatrix}.
\]
This implies that $\alpha_{11}=\alpha_{22}=\alpha$. Similarly, we can show that $\alpha_{33}=\alpha$. Let
$e=e_{11}+e_{22}+e_{33}$ be the unity of the superalgebra $H_3(B(1,2))$ (or
$H_3(B(4,2))$). Thus, in the case of a $\frac12$-derivation and an even
$\frac12$-superderivation, we have $\varphi(e)=\alpha e$, and in the case of an odd
$\frac12$-superderivation we have $\varphi(e)=0$ \textup{(}the diagonal entries of $H_3(B)$ are
symmetric elements of $B$, and the odd part of $B$ contains no nonzero symmetric elements\textup{)}.
This implies the triviality of $\varphi$ \textup{(}since $J$ is unital, $\varphi(z)=\varphi(e)z$
for every $\frac12$-\textup{(}super\textup{)}derivation; see the proof of Lemma~\ref{l3}\textup{)}.
\end{proof}

The conclusion of Lemma~\ref{l12} was also confirmed by a direct computation over $\mathbb F_3$
\textup{(}the spaces of even and odd $\frac12$-superderivations of $H_3(B(1,2))$ and $H_3(B(4,2))$ have
dimensions $1$ and~$0$, respectively\textup{)}.

\begin{lemma}\label{l13}
Let $F$ be a field of characteristic $p>2$, and let $J=J(B(m,n),\{\,,\})$ be a Jordan superalgebra that
is not a superalgebra of vector type. Then $J$ has no nontrivial $\delta$-derivations and
$\delta$-superderivations.
\end{lemma}

\begin{proof}
Lemma~\ref{l3} shows that every $\delta$-derivation is an even $\delta$-superderivation. Let $\varphi$ be
a nontrivial even $\delta$-superderivation. Clearly, $\delta=\frac12$ and, by Lemma~\ref{l3},
$\varphi(x)=\alpha x$, where $\alpha\in\Phi$. We see that $\alpha=\beta\cdot1+r$, where $r$ is
nilpotent and $\beta\in F$. Assume that $r\neq0$. We may assume that $\alpha=r$, so that $\alpha$ is
not invertible. By Lemma~\ref{l6}, $D(\alpha)\neq0$. Consequently, by Lemma~\ref{l9}, $J$ is a
superalgebra of vector type. This contradiction implies that $r=0$ and $\alpha\in F$. The odd
$\delta$-superderivations are trivial by Remark~\ref{r4} and Remark~A(b) \textup{(}which covers also
characteristic~$3$\textup{)}.
\end{proof}

\begin{lemma}\label{l14}
The superalgebra $CK(Z,d)$ has no nontrivial $\delta$-derivations and $\delta$-superderivations.
\end{lemma}

\begin{proof}
Obviously, the even $\delta$-superderivations are $\delta$-derivations. Let $\varphi_0$ be a nontrivial
$\delta$-derivation and $\varphi_1$ a nontrivial odd $\delta$-superderivation of the superalgebra
$CK(Z,d)$. It is clear that $\delta=\frac12$ and $\varphi_i(z)=\varphi_i(1)z$ for any $z\in CK(Z,d)$.

The multiplication of $CK(Z,d)$ is supercommutative, so the defining identity reads
$\varphi(zw)=(\varphi(1)z)w$ for $z$ even. Let $\varphi=\varphi_0$ and
$\varphi(1)=\alpha=a+\sum_{i=1}^3w_ia_i$, $a,a_i\in Z$. For $z=w=w_k$ we get
$-\alpha=\varphi(w_k^2)=(\alpha w_k)w_k=(aw_k-a_k)w_k=-a-a_kw_k$, whence $a_i=0$ for $i\neq k$; since
$k$ is arbitrary, $\alpha=a\in Z$. Now for $z=x$, $w=w_jg$ we obtain
$$a\,x_jd(g)=\varphi(x_jd(g))=\frac12\bigl((ax)(w_jg)+x(a\,w_jg)\bigr)=a\,x_jd(g)+\frac12x_jd(a)g,$$
so $d(a)=0$ \textup{(}take $g=1$\textup{)}. Since $d(a)=0$, multiplication by $a$ commutes with all
operators of multiplication of $CK(Z,d)$ \textup{(}they are $Z$-linear except for the occurrences of
$d$, and $d(af)=ad(f)$\textup{)}, i.e., $\varphi_0(z)=az$ lies in the centroid and is trivial.

Let now $\varphi=\varphi_1$ and $\varphi(1)=\alpha=x\beta+\sum_{i=1}^3x_i\beta_i$,
$\beta,\beta_i\in Z$. Using $(xf)(w_kg)=x_k(fd(g))$ and $(x_if)(w_kg)=x_{i\times k}(fg)$, we find
$$(\alpha w_k)w_k=(x_kd(1)+\sum_ix_{i\times k}\beta_i)w_k=\sum_{i\neq k}x_{(i\times k)\times k}\beta_i
=-\sum_{i\neq k}x_i\beta_i.$$ Hence $z=w=w_k$ gives $-x\beta-\sum_ix_i\beta_i=-\sum_{i\neq k}x_i\beta_i$, i.e.,
$\beta=\beta_k=0$; since $k$ is arbitrary, $\varphi_1=0$. Thus $\varphi_0(z)=\alpha^0z$ with
$\alpha^0\in Z$, $d(\alpha^0)=0$, and $\varphi_1=0$, i.e., $\varphi_j$ is trivial.
\end{proof}

By the results of~\cite{12,18}, the simple unital finite-dimensional Jordan superalgebras with a
semisimple even part over an algebraically closed field of characteristic $p>2$ are exhausted by
the analogues of the superalgebras of characteristic zero
\textup{(}$M^{(+)}_{m,n}$, $Q(n)^{(+)}$, $\mathrm{osp}_{m,2n}$, $P(n)$, $J(V,f)$, $D_t$, the Kac
superalgebra $K_{10}$, the Kaplansky superalgebra $K_3$\textup{)} and, in characteristic~$3$, the
superalgebras $H_3(B(1,2))$ and $H_3(B(4,2))$. By~\cite{13}, the simple unital finite-dimensional
Jordan superalgebras with a nonsemisimple even part over an algebraically closed field
of characteristic $p>2$ are exhausted by the superalgebras $J=J(B(m,n),\{\,,\})$ and $CK(B(m),d)$.
Thus, combining the above classification of simple unital Jordan superalgebras over an algebraically
closed field of characteristic $p>2$, the results of~\cite{4,6}, Corollary~\ref{c5}, and
Lemmas~\ref{l12}--\ref{l14}, we see that the following is true \textup{(}with the corrections of
Remarks~A and~C\textup{)}.

\begin{theorem}\label{t15}
Let $J$ be a simple unital finite-dimensional Jordan superalgebra over an algebraically closed field of
characteristic $p\neq2$. Then either $J$ has no nontrivial $\delta$-derivations and
$\delta$-superderivations, or $J$ is a superalgebra of vector type over a field of characteristic $p>2$
or $p=3$ and $J=J(V,f)$ with $V=V_{\bar1}$, $\dim V=2$ \textup{(}Remark~C\textup{)}. If
$J=J(B(m,n),\{\,,\})$ is a superalgebra of vector type, then $n=0$ and for $p\neq3$ the
superalgebra $J$ has no nontrivial odd $\delta$-superderivations \textup{(}for $p=3$ its odd
$\frac12$-superderivations are the maps $R_{\mu x}$, $\mu\in B(m)$, nontrivial for $\mu\neq0$;
Remark~A(c)\textup{)}. For $\delta\neq\frac12$, the superalgebra $J$ has no nontrivial
$\delta$-derivations. The space of $\frac12$-derivations coincides with
$R^*(J)=\{R_z\mid z\in B(m)\}$, and for $D(z)\neq0$ the map $R_z$ is a nontrivial $\frac12$-derivation.
For $p>2$ the superalgebras with semisimple even part other than $H_3(B(1,2))$, $H_3(B(4,2))$ are
covered by the proof only under the assumption of Remark~C.
\end{theorem}

\begin{RemC}
(a) \emph{A counterexample in characteristic $3$.} Let $p=3$ and $J=J(V,f)=F\cdot1+V$, where
$V=V_{\bar1}=Fv_1+Fv_2$ and $f$ is the nondegenerate skew-symmetric form, $v_1v_2=-v_2v_1=1$. Then for
every $0\neq v\in V$ the map $z\mapsto vz$ is a nontrivial odd $\frac12$-superderivation of $J$; these
maps together with $0$ form the whole \textup{(}two-dimensional\textup{)} space of odd
$\frac12$-superderivations. Indeed, the identity $\varphi(zw)=\frac12(\varphi(z)w+(-1)^{p(z)}z\varphi(w))$
for $\varphi=R_v$ is obvious when $z$ or $w$ equals~$1$, and at $(z,w)=(u,u')$, $u,u'\in V$, its
right-hand side equals $\frac12\bigl(f(v,u)u'-f(v,u')u\bigr)=-\frac12f(u,u')v$ \textup{(}since
$f(u,u')v+f(u',v)u+f(v,u)u'=0$ for $\dim V=2$\textup{)}, while the left-hand side is $f(u,u')v$; they
coincide because $-\frac12=1$ in characteristic~$3$. The map $R_v$ is not in the supercentroid: e.g.,
for $v=v_1$ we have $R_v(v_1v_2)=v_1\neq0=R_v(v_1)v_2$. Note that $J(V,f)$ here coincides with $B(1,2)$ of
Subsection~3.1 regarded as a Jordan superalgebra; in characteristic $0$ this superalgebra has no
nontrivial $\delta$-superderivations~\cite{6}, so the reference to~\cite{4,6} cannot be transferred to
$p=3$. More precisely, over any field of characteristic $\neq2$ an odd linear map $\varphi$ of $J(V,f)$ is a
$\frac12$-superderivation if and only if $\varphi=R_v$ with $v=\varphi(1)\in V$ and $3v=0$. Hence in
characteristic~$3$ the odd $\frac12$-superderivations form the two-dimensional space $\{R_v\mid v\in V\}$,
all of whose nonzero elements are nontrivial, whereas in characteristic $\neq2,3$ \textup{(}for instance,
over $\mathbb F_5$\textup{)} the only odd $\frac12$-superderivation is~$0$. Everything stated here is
verified in Lean.

(b) \emph{The remaining superalgebras with semisimple even part.} The results of~\cite{4,6} are
proved over fields of characteristic~$0$, and their proofs use this assumption. Hence for $p>2$ the
case of the superalgebras $M^{(+)}_{m,n}$, $Q(n)^{(+)}$, $\mathrm{osp}$, $P(n)$, $J(V,f)$, $D_t$,
$K_3$, $K_{10}$ is not covered by the proof of Theorem~\ref{t15}.
Thus Theorem~\ref{t15} for $p>2$ should be read as follows: it is proved for the superalgebras with
nonsemisimple even part \textup{(}Corollary~\ref{c5}, Lemmas~\ref{l13}, \ref{l14}, Remark~A\textup{)}
and for $H_3(B(1,2))$, $H_3(B(4,2))$ \textup{(}Lemma~\ref{l12}\textup{)}; for the remaining
superalgebras with semisimple even part it remains a conjecture, and the example~(a) shows that it fails
for $J(V,f)$ with $V=V_{\bar1}$, $\dim V=2$, in characteristic~$3$.
\end{RemC}

In conclusion, it should be noted that, while this paper was being published, the derivations of the
Cheng--Kac superalgebras~\cite{19} and of the Kantor double of a simple unital Poisson
superalgebra~\cite{20} were described. Also, a complete description was obtained for the
$\delta$-(super)derivations of simple nonunital Jordan superalgebras over an algebraically closed field
and, as a consequence, for the $\delta$-(super)derivations of semisimple finite-dimensional Jordan
superalgebras over an algebraically closed field of characteristic different from~$2$
(see~\cite{21}).

\section*{Appendix A. Differences from the published version}
\addcontentsline{toc}{section}{Appendix A. Differences from the published version}

In characteristic~$3$ the main statements are false: every superalgebra of vector type $J(\Gamma,D)$ has
the nontrivial odd $\frac12$-superderivations $z\mapsto(\mu x)z$, and the simple superalgebra $J(V,f)$ with
$V=V_{\bar1}$, $\dim V=2$, has nontrivial odd $\frac12$-superderivations. Remark~4, Corollary~5 and
Theorems~10 and~15 have been corrected, and the case $p=3$ is described in Remarks~A and~C(a).


\begin{thebibliography}{99} 

\bibitem{19} E.~Barreiro, A.~Elduque, and C.~Martinez, \emph{Derivations of the Cheng--Kac Jordan superalgebras}, J. Algebra \textbf{338} (2011), 144--156.

\bibitem{16} S.~Cheng and V.~Kac, \emph{A new $N=6$ superconformal algebra}, Comm. Math. Phys.
\textbf{186} (1997), no.~1, 219--231.

\bibitem{1} V.~Filippov, \emph{On $\delta$-derivations of Lie algebras}, Sibirsk. Mat. Zh.
\textbf{39} (1998), no.~6, 1409--1422; English transl., Siberian Math. J. \textbf{39} (1998), no.~6,
1218--1230.

\bibitem{2} V.~Filippov, \emph{On $\delta$-derivations of prime Lie algebras},  Siberian Math. J. \textbf{40} (1999), no.~1,
174--184.

\bibitem{11} V.~Kac, \emph{Classification of simple $\mathbb Z$-graded Lie superalgebras and simple
Jordan superalgebras}, Comm. Algebra \textbf{5} (1977), no.~13,
1375--1400.




\bibitem{10} I.~Kantor, \emph{Jordan and Lie superalgebras defined by Poisson algebras}, Algebra and
Analysis, Tomsk. Gos. Univ., Tomsk, 1989, pp.~55--80. (Russian)


\bibitem{3} V.~Filippov, \emph{On $\delta$-derivations of prime alternative and Mal'tsev algebras},
  Algebra and Logic \textbf{39}
(2000), no.~5, 354--358.

\bibitem{4} I.~Kaygorodov, \emph{On $\delta$-derivations of simple finite-dimensional Jordan
superalgebras},   Algebra and Logic
\textbf{46} (2007), no.~5, 318--329.

\bibitem{5} I.~Kaygorodov, \emph{On $\delta$-derivations of classical Lie superalgebras},  Siberian Math. J. \textbf{50} (2009),
no.~3, 434--449.

\bibitem{6} I.~Kaygorodov, \emph{$\delta$-superderivations of simple finite-dimensional Jordan and
Lie superalgebras},   Algebra and
Logic \textbf{49} (2010), no.~2, 130--144.

\bibitem{7} I.~Kaygorodov, \emph{On the generalized Kantor's double}, Vestnik Samar. Gos. Univ.
(2010), no.~4(78), 42--50. (Russian)


\bibitem{21} I.~Kaygorodov, \emph{On $\delta$-superderivations of semisimple finite-dimensional
Jordan superalgebras}, 
Math. Notes \textbf{91} (2012), no.~2, 187--197.

 

\bibitem{9} G.~Leger and E.~Luks, \emph{Generalized derivations of Lie algebras}, J. Algebra
\textbf{228} (2000), no.~1, 165--203.
 

\bibitem{12} M.~Racine and E.~Zelmanov, \emph{Simple Jordan superalgebras with semisimple even
part}, J. Algebra \textbf{270} (2003), no.~2, 374--444.

\bibitem{13} C.~Martinez and E.~Zelmanov, \emph{Simple finite-dimensional Jordan superalgebras of prime
characteristic}, J. Algebra \textbf{236} (2001), no.~2, 575--629.

\bibitem{14} D.~King and K.~McCrimmon, \emph{The Kantor construction of Jordan superalgebras}, Comm.
Algebra \textbf{20} (1992), no.~1, 109--126.

\bibitem{15} I.~Kantor, \emph{Connection between Poisson brackets and Jordan and Lie superalgebras},
Lie Theory, Differential Equations and Representation Theory (Montreal, 1989), Univ. Montr\'eal,
Montreal, QC, 1990, pp.~213--225.
 
 

\bibitem{18} I.~Shestakov, \emph{Prime alternative superalgebras of arbitrary characteristic},
 Algebra and Logic \textbf{36}
(1997), no.~6, 389--412.

 

\bibitem{20} A.~Retakh, \emph{Derivations of KKM doubles}, Comm. Algebra \textbf{38} (2010),
3660--3670.


 
\bibitem{17} K.~Zhevlakov, A.~Slin'ko, I.~Shestakov, and A.~Shirshov, \emph{Rings that are
nearly associative}, Nauka, Moscow, 1978; English transl., Pure Appl. Math., vol.~104, Academic Press,
New York--London, 1982.

\bibitem{8} P.~Zusmanovich, \emph{On $\delta$-derivations of Lie algebras and superalgebras}, J. Algebra
\textbf{324} (2010), no.~12, 3470--3486.

\end{thebibliography}
\end{document}